\documentclass[letterpaper,11pt]{article}
\usepackage{amsmath}
\usepackage{amssymb}
\usepackage{geometry}
 \usepackage{hyperref}
 \usepackage{amssymb}
 \usepackage{amsmath}
 \usepackage[latin1]{inputenc}
 \usepackage[T1]{fontenc}
\usepackage{times}
\usepackage{epsfig}
 \usepackage{dsfont}
 \usepackage{natbib}
 \usepackage{mathrsfs}
\date{}

\newtheorem{theorem}{Theorem}

\begin{document}

\title{Approximations for general bootstrap of empirical processes with an application to kernel-type density estimation}
\author{Salim BOUZEBDA\footnote{e-mail: salim.bouzebda@upmc.fr} \hbox{ }and Omar EL-DAKKAK\footnote{e-mail: omar.eldakkak@gmail.com} \\
Laboratoire de Statistique Th\'{e}orique et Appliqu\'{e}e (L.S.T.A.)\\
Universit\'{e} Paris VI}
\date{}
\maketitle

\begin{abstract}
\noindent The purpose of this note is to provide an approximation for the
generalized bootstrapped empirical process achieving the rate in
\cite{KMT}. The proof is based on much the same arguments used in
\cite{Lajos22000}. As a consequence, we establish an approximation
of the bootstrapped kernel-type density estimator.\\


\end{abstract}

\section{Introduction and Main Results}
Let $X_1,X_2,\ldots$ be a sequence of independent, identically
distributed [i.i.d.] random variables with common distribution
function $F(t)=P(X_1\leq t)$. The empirical distribution function of
$X_1,\ldots,X_n$ is
\begin{equation}
F_n(t) = \frac{1}{n}\sum_{i=1}^n \mathds{1}\{X_i\leq t\},\quad -\infty < t< \infty,
\end{equation}
where $\mathds{1}\{A\}$ stands for the indicator function of the
event $A$. Given the sample $X_1, \ldots, X_n,$ let $X_1^*, \ldots,
X_m^*,$ be conditionally independent random variables with common
distribution function $F_n.$ Let
\begin{equation}
F_{m,n}(t) = \frac{1}{m}\sum_{i=1}^m \mathds{1}\{X_i^*\leq
t\}, \quad -\infty < t < \infty,
\end{equation}
denote the \emph{classical} Efron (or multinomial) bootstrap (see,
e.g. \cite{Efron79} and \cite{EfronTibshirani1993} for more details).
Define the \emph{bootstrapped empirical process},
$\widehat{\alpha}_{m,n}$, by
\begin{equation}\label{bootsefron}
\alpha_{m,n}(t):=\sqrt{n}(F_{m,n}(t)-F_n(t)), 
\quad -\infty < t < \infty.
\end{equation}
Among many other things, \cite{BickelFreedman1981} established weak convergence of the process in (\ref{bootsefron}), which enabled
them to deduce the asymptotic validity of the bootstrap method in
forming confidence bounds for $F(\cdot)$. \cite{Shorack1982} gave a
simple proof of weak convergence of the process in
(\ref{bootsefron}) [see also \cite{ShorackWellner1986}, Section
23.1]. The Bickel and Freedman result for $\alpha_{m,n}$
has been subsequently generalized for empirical processes based on
observations in $\mathds{R}^d$, $d>1$ as well as in very general
sample spaces and for various set and function-indexed random
objects [see, for example \cite{beran1984}, \cite{BeranMillar1986},
\cite{Beran1987}, \cite{Gaenssler1987P}, \cite{Lohse1987}]. This
line of research found its ``final results'' in the work of
\cite{Gine1989,Gine1990} and \cite{CsorgoMason89}.

\noindent By now, the bootstrap is a widely used tool and,
therefore, the properties of $\alpha_{m,n}(t)$ are of great interest
in applied as well as in theoretical statistics. In fact, several
procedures can actually be described in terms of the empirical
process $\alpha_n(t)$, the limit distributions being functionals of
$B(F(t))$, where $B$ is a Brownian bridge. The fact that the limits
may depend on the unknown distribution $F(t)$ makes it important
that good approximations of these limiting distributions be found
and that is where the bootstrap proved to be a very effective tool.
There is a huge literature on the application of the bootstrap
methodology to nonparametric kernel density and regression
estimation, among other statistical procedures, and it is not the
purpose of this paper to survey this extensive literature. This
being said, it is worthwhile mentioning that the bootstrap as per
Efron's original formulation (see \cite{Efron79}) presents some
drawbacks. Namely, some observations may be used more than once
while others are not sampled at all. To overcome this difficulty, a
more general formulation of the bootstrap has been devised: the {\it
weighted} (or {\it smooth}) bootstrap, which has also been shown to
be computationally more efficient in several applications. For a
survey of further results on weighted bootstrap the reader is
referred to \cite{Bertail95}. Exactly as for Efron's bootstrap, the question
of rates of convergence is an important one (both in probability and
in statistics) and has occupied a great number of authors (see
\cite{csorgorevesz1981}, \cite{KMT}  \cite{Lajos22000} and the references therein).
\vskip7pt
\noindent In this note, we will consider a version of the
Mason-Newton bootstrap (see \cite{Mason1992}, and the references
therein). As will be clear, this approach to bootstrap is very
general and allows for a great deal of flexibility in applications.
 Let $(X_n)_{n\geq1}$ be a sequence of i.i.d. random variables defined on a probability space $(\Omega,\mathcal{A}, \mathbb{P}).$
 We extend $(\Omega,\mathcal{A}, \mathbb{P})$ to obtain a probability space $(\Omega^{(\pi)},\mathcal{A}^{(\pi)},P)$. The latter will carry the independent sequences $(X_n)_{n \geq 1}$ and $(Z_n)_{n \geq 1}$ (defined below) and will be considered rich enough as to allow the definition of another sequence $(B_n^*)$ of Brownian bridges, independent of all the preceding sequences. The possibility of such an extension is discussed in detail in literature; the reader is referred, e.g., to \cite{csorgorevesz1981}, \cite{KMT} and \cite{phillipberkes1977}. In the sequel, whenever an almost sure property is stated, it will be tacitly assumed that it holds with respect the the p.m. $P$ defined on the extended space.
\vskip7pt
\noindent Define a sequence $(Z_n)_{n\geq1}$ of i.i.d.
replic{\ae} of a {\it strictly positive} random variable $Z$ with
distribution function $G(\cdot)$, independent of the $X_n$'s. In the
sequel, the following assumptions on the $Z_n$'s will prevail:
\begin{enumerate}
\item[(A1)] $\mathds{E}(Z) = 1$; \quad $\mathds{E}(Z^2) = 2$ (or,
equivalently, $\mathbf{Var}(Z) = 1$).
\item[(A2)] There exists an
$\varepsilon > 0$, such that $$\mathds{E}(e^{tZ}) < \infty
~~\mbox{for all}~~ |t| \leq \varepsilon.$$
\end{enumerate}
For all $n \geq 1$, let $T_n = Z_1 + \cdots+ Z_n$ and define the random weights,
\begin{equation}\label{poids}
\mathscr{W}_{i;n}:= \frac{Z_i}{T_n}, \qquad i=1,\ldots,n.
\end{equation}
The quantity
\begin{equation}
F^*_n(t) = \sum_{i=1}^n \mathscr{W}_{i;n}\mathds{1}\{X_i\leq t\},~
\mbox{ for } ~ -\infty < t < \infty.
\end{equation}
will be called {\it generalized (or weighted) bootstrapped empirical
distribution function}. Analogously, recalling the empirical process
based on $X_1,\ldots,X_n$,
\begin{equation}
\alpha_n(t) = n^{1/2}(F_n(t) - F(t)),~-\infty < t < \infty,
\end{equation}
define the corresponding {\it generalized (or weighted) bootstrapped
empirical process} by
\begin{equation}
\alpha^*_n(t) = n^{1/2}(F^*_n(t) - F_n(t)),~-\infty < t < \infty.
\end{equation}

\noindent The system of weights defined in (\ref{poids}) appears in
\cite{Mason1992}, p.1617 where it is shown that it satisfies
assumptions ($\mathscr{W}_I$), ($\mathscr{W}_{II}$) and
($\mathscr{W}_{III}$) on p.1612 of the same reference, so that all the
results therein hold for the objects to be treated in this note.
In particular, weak convergence for the process $\alpha_n^*$ to a
Brownian bridge is proved. For more results concerning this version of the the weighted boostrapped empirical process, we refer the reader to \cite{deheuvels2008}.
Note that, as a special case of the system of weights we are considering, one can obtain the one used for Bayesian bootstrap (see \cite{rubin1981}).
\vskip7pt
\noindent In what follows, we obtain a KMT rate of
convergence for this process in sup norm. More precisely, we
consider deviations between the generalized bootstrapped empirical
process $\{\alpha_n^*(t): t\in\mathds{R}\}$ and a sequence of approximating Brownian
bridges $\{B^*_n(F(t)):t\in \mathds{R}\} $ on $\mathds{R}$. Our main
result goes as follows.

\begin{theorem}\label{1er theorem} Let assumptions {\rm (A1)} and {\rm (A2)} hold. Then, it is possible to
define a sequence of Brownian bridges $\{B_n^*(y): 0 \leq y\leq  1
\}$ such that, for all $\varepsilon,\eta>0$, there exists $N=N(\varepsilon,\eta)$, such that, for all $n \geq N$ and all $x>0,$
\begin{equation}\label{th}
P\left(\sup_{-\infty<t<\infty}\left|\alpha_n^*(t)-
B_n^*(F(t))\right|>3n^{-1/2}(K_1\log n+x)\right)\leq K_2 \exp \left(-\frac{K_3x}{(1+\varepsilon)^2}\right)+\eta,
\end{equation}
where $K_1$,  $K_2$ and $K_3$ are positive
universal constants.
\end{theorem}
The proof of Theorem \ref{1er theorem} is given in Section
\ref{thhh2}.

\vskip7pt

\noindent {\bf Remark 1} Theorem \ref{1er theorem} implies the
following approximation of the weighted bootstrap:
\begin{eqnarray}
 \sup_{-\infty<t<\infty}|\alpha_n^*(t)- B_n^*(F(t))|=O_P\left(\frac{\log
n}{n^{1/2}}\right).
\end{eqnarray}

\vskip7pt

\noindent {\bf Remark 2} Theorem \ref{1er theorem} turns out be useful in obtaining
confidence bands for the distribution function of the sample data.
We formalize this idea as follows: for $0<\alpha<1$, one has
\begin{eqnarray}
\lim_{n\rightarrow \infty}P\left(\sup_{-\infty<t<\infty}\sqrt{n}| F_n(t)-F(t)| \leq c(\alpha)\right)=P\left(\sup_{-\infty<t<\infty}|B(F(t))|\leq
c(\alpha), \right).
\end{eqnarray}
Note that for each fixed $t$, $B(F(t))$ is a zero-mean Gaussian
random variable with covariance structure
$$E(B(F(t))B(F(s)))=F(t\wedge s)-F(t)F(s)$$
where $t\wedge s:=\min(t,s)$. In practice, $c(\alpha)$ can, of
course, not be computed since the covariance structure of $B(F(t))$ depends on
the unknown cdf $F$. Instead, suppose $(Z_1^{(1)},\ldots
,Z_n^{(1)}),\ldots ,(Z_1^{(N)},\ldots ,Z_n^{(N)})$ are $N$
independent vectors of i.i.d. copies of $Z$, sampled independently
of the $X_i$'s. Define the random variables
\begin{equation}
\psi^j:=\sup_{-\infty<t<\infty}\left|\alpha_{n,j}^*(t) \right|,
\qquad j=1,\ldots ,N,
\end{equation}

\noindent where $\alpha_{n,j}^*$ denotes the generalized
bootstrapped empirical process constructed with the sample
$(Z_1^{(j)},\ldots ,Z_n^{(j)})$, $j=1,\ldots ,N$. Theorem \ref{1er theorem} accounts for the use of the smallest $z>0$ such that
\begin{equation*}
\frac{1}{N}\sum_{i=1}^N\mathds{1}\left\{\psi^j\leq z\right\}\geq
1-\alpha.
\end{equation*}
as an estimator of $c(\alpha).$
\vskip7pt

\noindent A direct consequence of Theorem \ref{1er theorem} and Theorem 1.5 in \cite{Lajos22000} is the following approximation for $\alpha_n^*(\cdot)$ based on a Kiefer process
\begin{theorem}\label{2 theorem}
There is a Kiefer process $\{K(t; x); 0\leq t\leq 1; 0\leq x\leq
\infty\}$ such that
\begin{equation}
\max_{ 1\leq k\leq n} \sup_{ -\infty<t<\infty} \left|
\sum_{i=1}^k(\mathscr{W}_{i;n}-1/n)\mathds{1}\{X_i\leq t\}- K(F(t),
k)\right| = O_P(n^{1/4}(\log n)^{1/2}).
\end{equation}
\end{theorem}
\section{An application to kernel density estimation}
Let $X_1,\ldots,X_n$ be independent random replic{\ae} of a random
variable $X\in\mathds{R}$ with distribution function $F(\cdot)$. We
assume that the distribution function $F(\cdot)$ has a density
$f(\cdot)$ (with respect  to the Lebesgue measure in $\mathds{R}$).
First of all, we introduce a kernel density estimator of $f(\cdot)$. To this
end, let $K(\cdot)$ be a measurable function
 fulfilling the following conditions
\begin{enumerate}
    \item[(K1)] $K(\cdot)$ is of bounded variation and compactly supported on
    $\mathds{R}$;
    \item[(K2)] $K\geq 0$ and  $\int K(u) du=1.$
\end{enumerate}
Now, define the Akaike-Parzen-Rosenblatt kernel density estimator of
$f(\cdot)$ (see \cite{akaike1954}, \cite{parzen1962} and
\cite{rosenblatt1956}) as follows: for all $x \in \mathds{R}$,
estimate $f(x)$ by
 \begin{equation}\label{A}
 f_{n,h_n}(x)=
 \frac{1}{nh_n}\sum_{i=1}^{n} K\left( \frac{x-X_i}{h_n} \right),
\end{equation}
where $\{h_n: n\ge 1\}$ is a sequence of positive constants
satisfying the conditions
$$h_n\downarrow 0 \quad \mbox{and} \quad
nh_n\uparrow\infty, \quad \mbox{as} \quad n\rightarrow \infty.$$
Secondly, we define the bootstrapped version of
$f_{n,h_n}(\cdot)$, by setting for all $h_n>0$ and $x\in
\mathds{R}$,
\begin{equation}\label{AB}
 f^*_{n,h_n}(x)=
 \frac{1}{h_n}\sum_{i=1}^{n}\mathscr{W}_{i;n} K \left( \frac{x-X_i}{h_n} \right),
\end{equation}
where $\mathscr{W}_{i;n}$ is defined in (\ref{poids}). We will provide an approximation
rate for the following process
\begin{equation}\label{pros}
 \gamma_n^*(x)=\sqrt{ nh_n^2}\left(f^*_{n,h_n}(x)-f_{n,h_n}(x)\right),
\quad -\infty<x<\infty.
\end{equation}

\noindent The following theorem, proved in the next Section, shows that a
single bootstrap suffices to obtain the desired approximation for
non-parametric kernel-type density estimators.
\begin{theorem}\label{thboots}
Let conditions {\rm (A1), (A2), (K1)} and {\rm (K2)}  prevail. Then we can
define Brownian bridges $\{B_n^*(y): 0 \leq y\leq 1 \}$ such that almost
surely along $X_1,X_2,\ldots ,$ as $n$ tends to infinity, we have
\begin{equation} \label{giuaniin}
\sup_{-\infty<x< \infty} \left |\gamma_n^*(x)- \int K \left( \frac{x-s}{h_n} \right) {\rm d}B_n^*(F(s)) \right |= O_P\left(\frac{\log n}{\sqrt{n}} \right).
\end{equation}
If, moreover, we suppose boundedness of the unknown density, f, i.e. if we suppose the existence of $M>0$ such that $\sup_{-\infty <x< \infty}f(x) \leq M,$ then, almost surely along $X_1,X_2, \ldots,$ as $n$ tends to infinity, 
\begin{equation} \label{giuanoon}
 \sup_{-\infty<x< \infty} \left |\gamma_n^*(x)-B_n^*(F(x))\int K(t) {\rm d}t \right |= O_P\left(\frac{\log n}{\sqrt{n}} + h_n \sqrt{\log h_n^{-1}}\right).
\end{equation}
\end{theorem}

\vskip7pt
\noindent {\bf Remark 3.} Under appropriate conditions, and using the same arguments rehearsed in the proof of Theorem \ref{thboots}, it is possible to obtain an approximation of a smoothed version of $F_n^*.$
\section{Proofs}\label{thhh2}

\noindent{\bf Proof of Theorem \ref{1er theorem}.} In the sequel, we
will write $\|\cdot\|$ to indicate
$\sup_{-\infty<t<+\infty}|\cdot|$. We have that
$$\|\alpha_n^*(t)-B_n^*(F(t))\|=\|\sqrt{n}(F_n^*(t)-F_n(t))-B_n^*(F(t))\|.$$
Now, it is easily seen that
\begin{equation}
\sqrt{n}(F_n^*(t)-F_n(t))=\left(\frac{n}{T_n}\right)\left[\frac{1}{\sqrt{n}}\left(\sum_{i=1}^n Z_i\mathds{1}\{X_i \leq t\}-F(t)T_n+(F(t)-F_n(t))T_n\right)\right],
\end{equation}
so that
$$\|\alpha_n^*(t)-B_n^*(F(t))\| \leq S_1(n)+S_2(n)+S_3(n),$$
where
\begin{equation}\label{esseuno}
S_1(n) := \left( \frac{n}{T_n} \right) \left\| \frac{1}{\sqrt{n}} \left(
\sum_{i=1}^n Z_i \mathds{1}\{X_i\leq t\}-T_n F(t)
\right)-B_n^*(F(t)) \right\|,
\end{equation}
where
\begin{equation}\label{essedue} S_2(n):=\left( \frac{n}{T_n} \right)
\left\|\frac{T_n}{\sqrt{n}}(F(t)-F_n(t)) \right\|,
\end{equation}
and where
\begin{equation}\label{essetre}
S_3(n):=\left|\frac{n}{T_n}-1\right| \left\| B_n^*(F(t)) \right\|.
\end{equation}
We start by dealing with the term $S_3(n)$. We will treat the cases
$x>Cn$ and $x \leq Cn$ ($C$ being a strictly positive constant) separately. Fix $x>Cn$ arbitrarily. Union bound gives for
all $n$,
$$P \left( S_3(n) \geq n^{-1/2} (x+c \log n)\right) \leq P\left(S_4(n) \geq \frac{x}{2 \sqrt{n}} \right)+P\left( \|B_n^*(F(t))\| \geq \frac{x}{2 \sqrt{n}}\right),$$
where

\begin{equation}\label{essequattro}
S_4(n) := \left( \frac{n}{T_n} \right) \left \| B_n^*(F(t)) \right
\|.
\end{equation}
\noindent 
Now, it is known that, for all $n \geq 1$ and all $x>n \geq 1$,
there exists a positive constant $c_1$, such that
\begin{equation}\label{normabienne}
P\left( \|B_n^*(F(t))\| \geq \frac{x}{2\sqrt{n}}\right) \leq c_1 \exp \left(- \frac{x^2}{4n} \right) \leq
\exp \left(-\frac{x}{4} \right).
\end{equation}
On the other hand, since strong law of large numbers gives
$$\left | \frac{n}{T_n} - 1 \right |\stackrel{a.s.}{\rightarrow} 0,$$
for all $\varepsilon,\eta >0,$ there exists $N_1=N_1(\varepsilon,\eta),$
such that, for all $n \geq N_1$,
\begin{equation}\label{epsilon_eta}
P \left( \left | \frac{n}{T_n} - 1 \right | \in
(0,\varepsilon) \right) \geq 1- \eta.
\end{equation}
Consequently, denoting the law of $ \frac{n}{T_n} $ by
$\mathcal{L}_{ \frac{n}{T_n} }$, independence of the $Z_n$'s from
the $B_n$'s gives
\begin{eqnarray}
\nonumber P\left(S_4(n) \geq \frac{x}{2 \sqrt{n}}\right)&=& P \left( S_4(n) \geq
\frac{x}{2\sqrt{n}}, \left | \frac{n}{T_n} - 1 \right | \in
(0,\varepsilon)
\right)\\
\nonumber && \qquad + \hbox{ }P \left( S_4(n) \geq \frac{x}{2\sqrt{n}}, \left
| \frac{n}{T_n} - 1 \right | \not\in (0,\varepsilon)
\right)\\
\nonumber & \leq & P \left( \frac{n}{T_n} \| B_n^*(F(t)) \| \geq
\frac{x}{2\sqrt{n}} \hbox{ }| \hbox{ } \left | \frac{n}{T_n} - 1 \right | \in
(0,\varepsilon) \right)\\
\nonumber && \qquad + \hbox{ }P \left( \left | \frac{n}{T_n} - 1
\right | \not\in (0,\varepsilon) \right)\\
\nonumber & \leq & \int_{1-\varepsilon}^{1+\varepsilon} P
\left( \| B_n^*(F(t)) \| > \frac{x}{2\sqrt{ny^2}} \hbox{ } | \hbox{
} \frac{n}{T_n}=y \right) \mathcal{L}_{\frac{n}{T_n}}(\hbox{d}y) +
\eta\\
\nonumber & \leq & P \left( \| B_n^*(F(t)) \| >  \frac{x}{2\sqrt{n(1+\varepsilon)^2}}
 \right)+ \eta \\
\label{venduuno} & \leq & c_1 \exp \left( -\frac{x}{4(1+\varepsilon)^2} \right)+ \eta,
\end{eqnarray}
where, in the last inequality, we have used (\ref{normabienne}).
Combining (\ref{normabienne}) and (\ref{venduuno}), we have that,
for all $\varepsilon, \eta >0,$ there exists
$N_1=N_1(\varepsilon,\eta),$ such that, for all $n \geq N_1,$
\begin{equation} \label{s3prima}
P \left( S_3(n) \geq n^{-1/2} (x+c \log n)\right) \leq (1+c_1)
\exp \left(-\frac{x}{4(1+\varepsilon)^2} \right)+\eta.
\end{equation}
Now we turn to the case $0<x \leq Cn.$ Again, by the union bound,
\begin{equation} \label{s3seconda}
P \left( S_3(n) \geq n^{-1/2} (x+c \log n)\right) \leq P \left(
\left| \frac{n}{T_n}-1 \right| > \sqrt{\frac{x}{n}} \right) + P \left( \|
B_n^*(F(t)) \| > \sqrt{x} \right).
\end{equation}
Again by (\ref{normabienne}), we have that for all $n,$
\begin{equation} \label{normab_n_2}
P \left( \| B_n^*(F(t)) \|>\sqrt{x} \right) \leq c_1 \exp (-x/2).
\end{equation}
On the other hand, by (\ref{epsilon_eta}), for all $\varepsilon,\eta>0$, there exists $N_1=N_1(\varepsilon,\eta)$ such that for all $n \geq N_1,$
\begin{eqnarray}
\nonumber P \left( \left| \frac{n}{T_n}-1 \right|>\sqrt{\frac{x}{n}}
\right)&=&P \left( \left| \frac{n}{T_n}-1 \right|>\sqrt{\frac{x}{n}},
\left|
\frac{n}{T_n}-1 \right| \in (0,\varepsilon) \right)\\
\nonumber&& \qquad \qquad + \hbox{ } P \left( \left| \frac{n}{T_n}-1
\right|>\sqrt{\frac{x}{n}}, \left|
\frac{n}{T_n}-1 \right| \not\in (0,\varepsilon) \right)\\
\nonumber & \leq & P \left( \left( \frac{n}{T_n} \right) \left|
\frac{T_n}{n}-1 \right| > \sqrt{\frac{x}{n}}, \left| \frac{n}{T_n}-1
\right| \in
(0,\varepsilon) \right)+ \eta\\
\label{mava} & \leq & P \left(  \left|
\frac{T_n}{n}-1 \right| > \sqrt{\frac{x}{n(1+\varepsilon)^2}} \right) + \eta
\end{eqnarray}
Use Theorem 2.6 in Petrov (1995) to find constants $c_2$ and $c_3$
such that
\begin{equation}\label{petrov}
P \left( \left| \frac{T_n}{n}-1 \right| >
\sqrt{ \frac{x}{n(1+\varepsilon)^2}} \right) \leq c_2 \exp\left( 
 \frac{-c_3 x}{(1+\varepsilon)^2}  \right).
\end{equation}
Combining (\ref{normab_n_2}), (\ref{mava}) and (\ref{petrov}), and
plugging in (\ref{s3seconda}), we deduce the existence of positive
universal constants $c_4$ and $c_5$ such that
\begin{equation}\label{s3ancora}
P \left( S_3(n) \geq n^{-1/2} (x+c \log n)\right) \leq c_4
\exp\left(  \frac{-c_5x}{(1+\varepsilon)^2}  \right)
+ \eta,
\end{equation}
so that one concludes, from (\ref{s3prima}) and (\ref{s3ancora}),
that for all $\varepsilon,\eta>0,$ there exists
$N=N(\varepsilon,\eta),$ such that, for all $n\geq N,$ and all $x>0$
\begin{equation} \label{s3fine}
P \left( S_3(n) \geq n^{-1/2} (x+c \log n)\right) \leq c_6
\exp\left(  \frac{-c_7x}{(1+\varepsilon)^2}  \right)
+ \eta,
\end{equation}
for some universal constants $c_6$ and $c_7.$

\noindent The proof is concluded once we show the existence of
universal positive constants $c_8$, $c_9$, $c_{10}$ and $c_{11}$
such that, for all $\varepsilon,\eta>0,$ there exists
$N_2=N_2(\varepsilon,\eta),$ and $N_3=N_3(\varepsilon,\eta)$ such that, for all $n\geq N_2,$ and all $x>0$
\begin{equation} \label{perche}
P \left( S_1(n) \geq n^{-1/2} (x+c \log n)\right) \leq c_8
\exp\left(  \frac{-c_9x}{(1+\varepsilon)^2}  \right)
+ \eta,
\end{equation}
and for all $n \geq N_3$ and all $x>0,$
\begin{equation} \label{percome}
P \left( S_2(n) \geq n^{-1/2} (x+c \log n)\right) \leq c_{10}
\exp\left(  \frac{-c_{11}x}{(1+\varepsilon)^2}
\right) + \eta.
\end{equation}
Since
\begin{equation}\nonumber
S_1(n)= \left( \frac{n}{T_n} \right) \left \| \frac{1}{\sqrt{n}} \left(
\sum_{i=1}^n Z_i \mathds{1}\{X_i\leq t\}-T_n F(t)
\right)-B_n^*(F(t)) \right \|,
\end{equation}
formula (3.7) in \cite{Lajos22000} combined with arguments similar to those used for the
term $S_3(n)$ imply (\ref{perche}). As for (\ref{percome}), formula (3.5) in \cite{Lajos22000} together with the by now usual $\varepsilon,\eta$ argument conclude the proof. \hfill$\Box$\\
\vskip7pt
\noindent{\bf Proof of Theorem \ref{thboots}.} We start by proving (\ref{giuaniin}). We have for $x\in \mathds{R}$
\begin{eqnarray}
\sqrt{nh_n^2}\left(\widehat{f}^*_{n,h_n}(x)-\widehat{f}_{n,h_n}(x)\right)&=&\int
K\left((x-
 s)/h_n\right) \hbox{d} \{n^{1/2}(F^*_n(s)-F_n(s))\}\nonumber\\&=& \int K\left((x-
 s)/h_n\right)\hbox{d}\alpha_n^*(s)\nonumber.
\end{eqnarray}
Integration by parts implies that
\begin{eqnarray}\label{17}
\int K\left(\frac{x-s}{h_n}\right)\hbox{d}\alpha_n^*(s)=-\int \alpha_n^*(x-th_n)\hbox{d}K(t),
\end{eqnarray}
and
\begin{equation} \label{17bis}
\int K\left(\frac{x-s}{h_n}\right)\hbox{d} B_n^*(F(s))=-\int B_n^*(F(x-th_n))\hbox{d}K(t). 
\end{equation}
Now, Theorem \ref{1er theorem} together with condition (K1) give
\begin{eqnarray}
\nonumber\lefteqn{\sup_{-\infty<x<\infty}\left|\int \alpha_n^*(x-th_n)\hbox{d}K\left(t\right)-\int B_n^*(F(x-th_n))\hbox{d}K\left(t\right)\right|}\\
&& \qquad \qquad \qquad\leq \sup_{-\infty<u<\infty}|\alpha_n^*(u)-B_n^*(F(u))|\int
\hbox{d}|K\left(t\right)| =O_P\left(\frac{\log n}{\sqrt{n}}\right),
\end{eqnarray}
thus proving (\ref{giuaniin}).

\noindent Once (\ref{giuaniin}) is at hand, to prove (\ref{giuanoon}), it suffices to bound
\begin{eqnarray}\label{18}
\left|\int B_n^*(F(x-th_n))\hbox{d}K\left(t\right)-B_n^*(F(x))\right|\leq
\int \left| B_n^*(F(x-th_n))-B_n^*(F(x))\right|\hbox{d}K\left(t\right),
\end{eqnarray}
in probability. By condition (K1), and provided the unknown density $f$ is bounded (by a strictly positive constant, say $M$), for $n$ large enough,
\begin{eqnarray}\label{19}
\left| B_n^*(F(x-th_n))-B_n^*(F(x))\right|\leq \sup_{|u-v|\leq
\delta_n}\left| B_n^*(u)-B_n^*(v)\right|
\end{eqnarray}
where $\delta_n=Mh_n$. Now, it is always possible to define a Brownian Bridge, $\{B^*(y):0 \leq y\leq  1\}$, on the
same probability space carrying the sequence of Brownian Bridges $\{B_n^*(y):0 \leq y\leq  1\}_{n \geq 1}$, such that for all $n,$ and all $\varepsilon>0$
\begin{eqnarray*}
\lefteqn{P\left(\{2\delta_n\log \delta_n^{-1}\}^{-1/2}\sup_{|u-v|<h}\sup_{h\in[0,\delta_n]}|B_n^*(u)-B_n^*(v)|>1+\varepsilon\right)}\\
&=&P\left(\{2\delta_n\log
\delta_n^{-1}\}^{-1/2}\sup_{|u-v|<h}\sup_{h\in[0,\delta_n]}|B^*(u)-B^*(v)|>1+\varepsilon\right).
\end{eqnarray*}
Since $\delta_n\rightarrow 0$, by Theorem 1.4.1 in
\cite{csorgorevesz1981}, we have with probability one
\begin{equation}
\lim_{n\rightarrow\infty}\{2\delta_n\log
\delta_n^{-1}\}^{-1/2}\sup_{|u-v|<h}\sup_{h\in[0,\delta_n]}|B^*(u)-B^*(v)|=1.
\end{equation}
Thus, as $n\rightarrow \infty$,
\begin{equation*}P\left(\{2\delta_n\log \delta_n^{-1}\}^{-1/2}\sup_{|u-v|<h}\sup_{h\in[0,\delta_n]}|B_n^*(u)-B_n^*(v)|>1+\varepsilon\right)\rightarrow0,
\end{equation*}
giving
\begin{equation}\label{21}
\sup_{|u-v|\leq
h}\sup_{h\in[0,\delta_n]}|B^*_n(u)-B^*_n(v)|=O_P \left(\sqrt{2\delta_n\log
\delta_n^{-1}} \right).
\end{equation}
Put (\ref{17}), (\ref{17bis}), (\ref{18}), (\ref{19}) and (\ref{21}) together to obtain
\begin{equation*}
 \sup_{-\infty<x< \infty} \left |\gamma_n^*(x)-B_n^*(F(x))\int {\rm d}K(t)\right |= O_P\left(\frac{\log n}{\sqrt{n}}+ h_n \sqrt{\log h_n^{-1}} \right),
\end{equation*}
thus completing the proof of Theorem. \hfill$\Box$

\vskip7pt
\noindent {\bf Acknowledgement} The authors are grateful to Professors 
Paul Deheuvels and Giovanni Peccati for a number of very useful discussions.

\end{document}